\documentclass[10pt]{amsart}
\usepackage{amsfonts,amssymb,amscd,amsmath,enumerate,verbatim,calc}
\usepackage{color,soul}
\usepackage{mathtools}
\usepackage{graphicx}
\everymath{\displaystyle}
\usepackage{fancyhdr}
\usepackage{tikz}
\usepackage{url}
\usetikzlibrary{shapes.geometric, arrows}
\tikzstyle{arrow} = [thick,->,>=stealth]
\tikzstyle{process} = [rectangle, minimum width=4cm, minimum height=2cm, text centered, text width=2.5cm, draw=green, fill=green!10]
\theoremstyle{plain}

\theoremstyle{definition}

\begin{document}
\title{Almost $SC^*$-normal spaces}
\author{Neeraj Kumar Tomar}
\email{neer8393@gmail.com}
\address{Department of Applied Mathematics, Gautam Buddha University, Greater Noida, Uttar Pradesh 201312, India}

\author{Fahed Zulfeqarr}
\email{fahed@gbu.ac.in}
\address{Department of Applied Mathematics, Gautam Buddha University, Greater Noida, Uttar Pradesh 201312, India}

\author{Amit Ujlayan}
\email{amitujlayan@gbu.ac.in}
\address{Department of Applied Mathematics, Gautam Buddha University, Greater Noida, Uttar Pradesh 201312, India}

\date{\today}
	
\begin{abstract}
This paper introduces and investigates a new class of almost normal spaces, referred to as almost $SC^*$-normal spaces, which are defined using $SC^*$-open sets. Building on the work of A. Chandrakala and K. Bala Deepa Arasi, we explore several properties of these spaces within the framework of topology. Moreover, we obtain some new characterizations and preservation theorems of almost $SC^*$-normal spaces.

\end{abstract}
\maketitle

\section*{Introduction} 
In this paper, we introduced the concept of almost $SC^*$-normal by using $SC^*$-open set due to A. Chandrakala and K. Bala Deepa Arasi [9] and obtained several properties of such a space. In 1970, Singal and Arya\cite{singal1970almost} introduced the concept of almost normal spaces as a generalization of normal - spaces by using regularly closed sets and obtained several properties of such a space. Recently, Hamant Kumar and M. C. Sharma\cite{Hamantalmost2022} introduced a new class of spaces, namely almost $\gamma$-normal and mildly $\gamma$-normal spaces are weaker form of $\gamma$- normal spaces. In this paper, we explore the properties of almost $\gamma$-normal and mildly $\gamma$-normal spaces, demonstrating that they are regularly open hereditary. We explore the connections between these spaces and certain weaker variants of almost normal and mildly normal spaces. 

Furthermore, we introduce and analyze the notions of of $gSC^*$-closed sets, $rgSC^*$-closed sets, regularly $SC^*$-closed sets, $T$-$SC^*$-closed sets, $T$-$SC^*$-open sets, and almost $SC^*$-irresolute functions. Additionally, we establish new characterizations and preservation theorems for various topological properties, including almost $SC^*$-normal spaces by exploring their structural relationships and interactions with other generalized forms of normality. Throughout this work, $(X, \tau)$ and $(Y, \sigma)$ denote topological spaces, with no specific separation axioms assumed unless explicitly stated.

\section{Preliminaries and Notations} 
In the following discussion, all spaces refer to topological spaces without any assumed separation axioms unless explicitly stated. A function $f :(X,\tau)\rightarrow (Y,\sigma)$ (or simply $f : X\rightarrow Y $), represents a mapping between two topological spaces $(X,\tau)$ and $(Y,\sigma)$. For any subset $A$ of a topological space $X$. The closure and the interior of $A$ are denoted by $\overline{A}$ and $A^\circ$, respectively.

\subsection{Definition:}
A subset $A$ of a topological space $X$ is said to be 

$(1)$ regular open ($r$-open) \cite{stone1937applications} if $A=int(cl(A))$.

$(2)$ semi open ($s$-open) \cite{levine1963semi} if $A\subset cl(int(A))$.

$(3)$ $\alpha$-open \cite{njȧstad1965some} if $A\subset (int(cl(int(A)))$.

$(4)$ $g\alpha$-closed\cite{maki1993generalized} if $\alpha $-$cl(A)\subset U$, whenever $A\subset U$, and $U$ is $\alpha$-open in $X$.

$(5)$ $rg\alpha$-closed\cite{vadivel2009rgalpha} if $\alpha $-$cl(A)\subset U$, whenever $A\subset U$, and $U$ is regularly $\alpha$-open in $X$.

$(6)$ $c^*$-open \cite{malathi2017pre} if $int(cl(A))\subset A\subset cl(int(A))$.

The complement of a $r$-open set is referred to as a $r$-closed set, while the complements of $s$-open and $c^*$-open sets are called $s$-closed and 
$c^*$-closed sets, respectively.

 For any subset $A$ of a topological space, the $c^*$-closure (resp. $s$-closure, $\alpha$-closure, $g\alpha$-closure, $rg\alpha$-closure) of $A$ is defined as the intersection of all $c^*$-closed (resp. $s$-closed, $\alpha$-closed, $g\alpha$-closed, $rg\alpha$-closed) sets that contain $A$ and is denoted by $c^*$-$\overline{A}$ (resp. $s$-$\overline{A}$, $\alpha$-$\overline{A}$, $g\alpha$-$\overline{A}$, $rg\alpha$-$\overline{A}$. Similarly, the $c^*$-interior (resp. $s$-interior, $\alpha$-interior, $g\alpha$-interior, $rg\alpha$-interior) of $A$ is defined as the union of all $c^*$-open (resp. $s$-open, $\alpha$-open, $g\alpha$-open, $rg\alpha$-open) sets contained in $A$. This is denoted by $c^*$-$A^\circ$ (resp. $s$-$A^\circ$, $\alpha$-$A^\circ$, $g\alpha$-$A^\circ$, $rg\alpha$-$A^\circ$).\\

\subsection{Definition} A subset $A$ of a space $X$ is said to be $SC^*$-closed\cite{Chandrakala2022} if $scl(A)\subset U$ whenever $A\subset U$ and $U$ is $c^*$-open in $X$. The complement of $SC^*$-closed set is said to be $SC^*$-open.

In a topological space, the $SC^*$-closure of a subset $A$ is the smallest $SC^*$-closed set that includes $A$. It is represented as $SC^*$-$\overline{A}$ and is obtained by taking the intersection of all $SC^*$-closed sets that contain $A$. Likewise, the $SC^*$-interior of a subset $A$, denoted as $SC^*$-$A^\circ$, is defined as the largest $SC^*$-open set contained within $A$. It is obtained by taking the union of all $SC^*$-open sets that are subsets of $A$.\\
The collection of all $SC^*$-open sets in a topological space $X$ is represented by $SC^*O(X)$, while the families of $SC^*$-closed, $r$-open, $r$-closed, $s$-open, and $s$-closed sets are denoted by $SC^*C(X)$ (resp. $R$-$O(X)$, $R$-$C(X)$, $S$-$O(X)$, $S$-$C(X)$). 

\subsection{Definition} A subset $A$ of a topological space $X$ is said to be 
\begin{flushleft}
    
$(1)$ {$g$-closed} \cite{levine1970gen} if $cl(A)\subset U$ whenever $A\subset U$ and 
$U\in \tau$. 

$(2)$ {generalized $SC^*$-closed}\cite{tomar2024scnormal} (briefly $gSC^*$-closed)  if  $SC^*$-$cl(A)\subset U$ whenever $A\subset U$ and $U\in\tau$.

$(3)$ {$SC^*$generalized-closed}\cite{tomar2024scnormal} (briefly $SC^*g$-closed)  if  $SC^*$-$cl(A)\subset U$ whenever $A\subset U$ and $U\in SC^*O(X)$.

$(4)$ regularly $SC^*$-open if there is a regularly open set $U$ such that $U\subset A\subset SC^*$-$cl(U)$.

$(5)$ $rgSC^*$-closed if $SC^*$-$cl(A)\subset U$, whenever $A\subset U$, and $U$ is regularly$SC^*$-open in $X$.\\

The complement of a $g$-closed set in a topological space is referred to as a $g$-open set. Similarly, the complements of $gSC^*$-closed, $SC^*g$-closed, and $rgSC^*$-closed sets are called $gSC^*$-open, $SC^*g$-open, and $rgSC^*$-open sets, respectively.

The complement of regularly $SC^*$-open set to be regularly $SC^*$-closed set.
\end{flushleft}

\subsection{Remark.} The following relationships hold for different types of subset properties in a topological space:

\vspace{8.0mm}

\begin{flushleft}
    closed $\Rightarrow$ $SC^*$-closed  $\Leftrightarrow$ $SC^*g$-closed  $\Leftrightarrow$  $gSC^*$- closed \\$\Downarrow$\\ $g$-closed
\end{flushleft}
\vspace{.4cm}
This representation outlines the logical implications between these closure properties.

\vspace{7.5mm}

\subsection{Remark.} The following relationships hold for different types of subset properties in a topological space:

\vspace{8.0mm}

\begin{flushleft}
   closed $\Longrightarrow$  $\alpha$-closed $\Longrightarrow$ $g\alpha$-closed  $\longrightarrow$  $rg\alpha$- closed\\ \hspace{2.2cm}$\Downarrow\hspace{1.8cm}\Downarrow\hspace{2cm}\Downarrow$\\ \hspace{1.6cm} $SC^*$-closed\hspace{.4cm} $gSC^*$-closed\hspace{.5cm}$rgSC^*$-closed
\end{flushleft}
\vspace{.4cm}
This representation outlines the logical implications between these closure properties.
\vspace{7.5mm}

\subsection{Example} Consider the set $X$ = $\{i,j,k,l \}$ with the topology $\tau = \{ \phi, \{i\}, \{j\}, \{i,j\},\{i,j,k\},X \}$. The subset $A = \{k\}$ is $\alpha$-closed and $SC^*$-closed but may not necessarily be closed in the topological space $X$ .

\subsection{Example} Consider $X$ = $\{i,j,k,l \}$ with the topology $\tau = \{ \phi, \{i\}, \{j\}, \{i,j\},\{i,j,k\},X \}$. The set $A = \{k\}$ is both $\alpha$-closed and $gSC^*$-closed set but it is not a closed set in $X$.

\subsection{Example} Consider the set $X$ = $\{i,j,k,l\}$ with the topology $\tau = \{\phi,\{j,l\},\{i,j,l\},\{j,k,l\},X\}$. Then

\begin{flushleft}
    $(1)$ closed sets in $(X,\tau)$ are $\phi$, $X$, $\{i\}$, $\{k\}$, $\{i,k\}$.

 $(2)$ $g$-closed set in $(X,\tau)$ are $\phi$, $X$, $\{i\}$, $\{k\}$, $\{i,k\}$, $\{i,j,k\}$, $\{i,k,l\}$.

 $(3)$ $SC^*$-closed set in $(X,\tau)$ are $\phi$, $X$, $\{i\}$, $\{j\}$, $\{k\}$, $\{l\}$, $\{i,j\}$, $\{i,k\}$, $\{i,l\}$, $\{j,k\}$, $\{j,l\}$, $\{k,l\}$, $\{i,j,k\}$, $\{i,j,l\}$, $\{i,k,l\}$, $\{j,k,l\}$.

$(4)$ $gSC^*$-closed set in $(X,\tau)$ are $\phi$, $X$, $\{i\}$, $\{j\}$, $\{k\}$, $\{l\}$, $\{i,j\}$, $\{i,k\}$, $\{i,l\}$, $\{j,k\}$, $\{j,l\}$, $\{k,l\}$, $\{i,j,k\}$, $\{i,j,l\}$, $\{i,k,l\}$, $\{j,k,l\}$.

$(5)$ $SC^*g$-closed set in $(X,\tau)$ are $\phi$, $X$, $\{i\}$, $\{j\}$, $\{k\}$, $\{l\}$, $\{i,j\}$, $\{i,k\}$, $\{i,l\}$, $\{j,k\}$, $\{j,l\}$, $\{k,l\}$, $\{i,j,k\}$, $\{i,j,l\}$, $\{i,k,l\}$, $\{j,k,l\}$.
\end{flushleft}

\section{almost $SC^*$-normal spaces}
\subsection{Definition} A topological space $X$ is said to be almost normal\cite{singal1970almost} (resp. almost $SC^*$-normal) if for every pair of disjoint closed sets $A$ and $B$, one of which is closed and other is regularly closed, there exist disjoint open(resp.SC*-open)sets $U$ and $V$ of $X$ such that $A\subset U$ and $B\subset V$.

\subsection{Example} Let $X = \{i,j,k,l\}$ and $\tau = \{\phi, \{i\}, \{j\}, \{i,j\}, \{k,l\}, \{i,k,l\}, \{j,k,l\}, X \}$. In this space, the subset $A = \{j\}$ is closed, while the subset $B = \{i\}$ is $r$-closed. There are separate open sets $U = \{j,k,l\}$ and $V = \{i\}$ where $A$ is contained in $U$ and $B$ is contained in $V$, ensuring they remain disjoint. Thus, $X$ satisfies the conditions for being almost normal. Moreover, since every open set in $X$ is also $SC^*$-open, the space is almost $SC^*$-normal.

A normal space implies an almost normal space, which in turn implies an almost $SC^*$-normal space. This relationship can be represented as follows:

\vspace{3.7mm}

normal $\Rightarrow$  almost-normal  $\Rightarrow$ almost $SC^*$-normal

\vspace{3.9mm}

\subsection{Lemma} If $A$ be a subset of a topological space $X$ is $rgSC^*$-open if and only if for every regularly closed set $F$ with $F\subset A$,  it follows that $F\subset SC^*$-$A^\circ$.

 \subsection{Theorem}The following conditions are equivalent for any topological space $X$:
\begin{flushleft}
   
$(1)$ $X$ is almost $SC^*$-normal,

$(2)$ If $I$ is a closed set and $J$ is a regularly closed set, then there exist disjoint $gSC^*$-open sets $M$ and $N$ such that $I\subset M$ and $J\subset N$.

$(3)$ If $I$ is a closed set and $J$ is a regularly closed set, then there exist disjoint $rgSC^*$-open sets $M$ and $N$ such that $I\subset M$ and $J\subset N$.

$(4)$ If $I$ is a closed set and $J$ is a regularly open set containing $I$, then there exist a $gSC^*$-open set $M$ satisfying $I\subset M\subset SC^*$-$\overline{I}\subset J$.

$(5)$ If $I$ is a closed set and $J$ is a regularly open set containing $I$, then there exist a $rgSC^*$-open set $M$ satisfying $I\subset U\subset SC^*$-$\overline{I}\subset J$.

$(6)$ If $I$ and $J$ are disjoint sets, where one is closed and the other is regularly closed, then there exist disjoint  $SC^*$-open sets $M$ and $N$ such that $I\subset M$ and $J\subset N$ and $M\cap N$ = $\phi $

\textbf{Proof} We will prove the implications step by step.\\$(1)\Rightarrow (2)$, $(2)\Rightarrow (3)$, $(4)\Rightarrow (5)$, $(3)\Rightarrow (4)$, $(5)\Rightarrow (6)$, $(6)\Rightarrow (1)$.\\
$(1)\Rightarrow (2)$. Assume that $X$ is almost $SC^*$-normal. Let $I$ be a closed set and $J$ be a regularly closed set in $X$. By the definition of almost $SC^*$-normality, there exist disjoint $SC^*$-open sets $M$ and $N$ such that $I\subset M$ and $J\subset N$. Since every $SC^*$-open set is also $gSC^*$-open, it follows that $M$ and $N$ are $gSC^*$-open sets, satisfying the required condition $I\subset M$ and $J\subset N$.

\vspace{1.2mm}
$(2)\Rightarrow (3)$. Let $I$ be a closed set and $J$ be a regularly closed set in $X$. By the assumption in \textbf{(2)}, there exist disjoint $gSC^*$-open sets $M$ and $N$ such that $I\subset M$ and $J\subset N$.
Since every $gSC^*$-open set is also an $rgSC^*$-open set, it follows that $M$ and $N$ are $rgSC^*$-open sets. Thus, we still have $I\subset M$ and $J\subset N$, satisfying the condition in \textbf{(3)}.

\vspace{1.2mm}
$(4)\Rightarrow (5)$. Let $I$ be a closed set and $J$ be a regularly open set such that $I\subset J$. By the assumption in \textbf{(4)}, there exists a $gSC^*$-open set $M$ in $N$ satisfying: $I\subset M\subset SC^*$-$\overline{I}\subset J$
Since every $gSC^*$-open set is also an $rgSC^*$-open set, the set $M$ retains the same containment properties while being $rgSC^*$-open. Thus, there exists an $rgSC^*$-open set $M$ such that: $I\subset M\subset SC^*$-$\overline{I}\subset J$. This confirms the validity of \textbf{(5)}.

\vspace{1.2mm}
$(3)\Rightarrow (4)$. Let $I$ be a closed set and $J$ be a regularly open set such that $I\subset J$. By the assumption in \textbf{(3)}, there exist disjoint $rgSC^*$-open sets $M$ and $N$ such that: $I\subset M$ and $X - J\subset N$. From \textbf{Lemma 2.3}, we know that: $X - J\subset SC^*$-$N^\circ$.
Additionally, the property $SC^*$-$\overline{M} \cap SC^*$-$N^\circ$ = $\phi$ holds. Using these facts, we conclude: $I\subset M\subset SC^*$-$\overline{M}\subset X - SC^*$-$N^\circ\subset J$
Thus, the statement in \textbf{(4)} follows.

\vspace{1.2mm}
$(5)\Rightarrow (6)$. Consider a closed set $I$ and a regularly open set $J$ that contains $I$. This implies that the complement $X - J$ is a regularly closed set.
By the assumption in \textbf{(5)}, there exists an $rgSC^*$-open set $K$ in $X$ such that: $I\subset K\subset SC^*$-$\overline{K}\subset X - J$.
Define: $M$ = $SC^*$-$K^\circ$,\\ $N$ = $X - SC^*$-$\overline{K}$, Since $M$ and $N$ are disjoint $SC^*$-open sets in $X$, it follows that $I\subset M$ and $J\subset N$, as required.

\vspace{1.2mm}
$(6)\Rightarrow (1)$. This follows directly from the definition of almost $SC^*$-normality, making the implication straightforward.

\end{flushleft}

\subsection{Definition} A function $f: X\rightarrow Y $ is called \textbf{$rc$-continuous}[10] if for each regular closed set $F$ in $Y$, $f^{-1}(F)$ is regularly closed in $X$.
\subsection{Definition} A function $f: X\rightarrow Y $ is termed \textbf{$T$-$SC^*$-open} (resp. \textbf{$T$-$SC^*$-$closed$}) if the image of every $SC^*$-open (resp. $SC^*$-closed) set in $X$ remains $SC^*$-open (resp. $SC^*$-closed) in $Y$.

\subsection{Definition} A function $f: X\rightarrow Y $ is said to be  \textbf{almost $SC^*$-irresolute} if, for every point $x\in X$ and for any $SC^*$-neighborhood $N$ of $f(x)$, the $SC^*$-closure of the preimage $f^{-1}(N)$ serves as an $SC^*$-neighborhood of $x$.

\section{Preservation theorems}

\subsection{Theorem} Let $f: X\rightarrow Y $ be a surjective function that is continuous, $T$-$SC^*$-open, $rc$-continuous and almost $SC^*$-irresolute. If $X$ is an almost $SC^*$-normal space, then $Y$ is almost $SC^*$-normal.

\textbf{Proof} Let $I$ be a closed set in $Y$ and $J$ be a regularly open set such that $I\subset J$. Since 
$f$ is $rc$-continuous, the preimage $f^{-1}(I)$ is a closed set in $X$, and it is contained in the regularly open set $f^{-1}(J)$. Since $X$ is almost $SC^*$-normal, \textbf{Theorem 2.4}, ensures the existence of an $SC^*$-open set $N$ in $X$ such that: $f^{-1}(I)\subset N\subset SC^*$-$cl(N)\subset f^{-1}(J)$ Applying $f$ to both sides, we obtain:
$f(f^{-1}(I))\subset f(N)\subset f(SC^*$-$cl(N))\subset f( f^{-1}(J))$. Since $f$ is an $T$-$SC^*$-open and almost $SC^*$-irresolute surjection, it follows that $f(N)\in SC^*O(Y)$, Thus, we conclude:\\ $I\subset f(N)\subset SC^*$-$cl(f(N))\subset J$. By \textbf{Theorem 2.4}, this confirms that $Y$ is almost $SC^*$- normal, completing the proof.

\subsection{Theorem} If $f:X\rightarrow Y$ is an $rc$-continuous and $T$-$SC^*$-closed mapping  from an almost $SC^*$-normal space $X$ onto $Y$, then $Y$ is also almost $SC^*$-normal.

\textbf{Proof} The proof follows directly from the definitions and properties of $rc$-continuity, $T$-$SC^*$-closed maps, and almost $SC^*$-normal spaces. By applying these properties, it is straightforward to verify that $Y$ inherits the almost $SC^*$-normality from $X$.


\bibliographystyle{amsplain}
\bibliography{references}

\end{document}